\documentclass[12pt]{amsart}
\usepackage{amssymb,times,fullpage}

\newtheorem{thm}{Theorem}
\newtheorem{prop}{Proposition}
\newtheorem{lem}{Lemma}

\newtheorem{conj}{Conjecture}

\theoremstyle{remark}
\newtheorem{rem}{Remark}

\theoremstyle{definition}
\newtheorem{defn}{Definition}

\newcommand{\C}{\mathbb{ C}}

\newcommand{\R}{\mathbb{ R}}

\newcommand{\Z}{\mathbb{ Z}}

\newcommand{\Spc}{\mathrm{Spin}^c}

\newcommand{\dvol}{{d\text{vol}}}

\newcommand{\FF}{\mathcal F}
\newcommand{\PP}{\overline{\lambda}}

\newcommand{\fs}{\mathfrak s}

\newcommand{\DA}{{D_{\! A}^+}}

\title{Monopole classes and Perelman's invariant of four-manifolds}
\author{D.~Kotschick}
\address{Mathematisches Institut, Ludwig-Maximilians-Universit\"at M\"unchen,
Theresienstr.~39, 80333 M\"unchen, Germany}
\email{dieter@member.ams.org}
\date{August 20, 2006; MSC 2000 classification: primary 57R57; secondary 53C25, 57R55}

\begin{document}

\begin{abstract}
We calculate Perelman's invariant for compact complex surfaces and a few other
smooth four-manifolds. We also prove some results concerning the dependence of 
Perelman's invariant on the smooth structure.
\end{abstract}

\maketitle

\section{Introduction}

In his celebrated work on the Ricci flow~\cite{P1,P2}, G.~Perelman introduced an interesting invariant of closed manifolds
of arbitrary dimension. By definition, Perelman's invariant is closely related to the Yamabe invariant or sigma constant. 
For three-manifolds that do not admit a metric of positive scalar curvature, Perelman's work~\cite{P1,P2} shows that his 
invariant is equivalent to M.~Gromov's minimal volume~\cite{gromov}, which is a priori a very different kind of invariant. 

In this paper we calculate Perelman's invariant for compact complex surfaces and show that it essentially coincides with the 
Yamabe invariant: 
\begin{thm}\label{t:cx}
Let $Z$ be a minimal complex surface with $b_1(Z)$ even and with $b_2^+(Z)>1$. Then Perelman's
invariant for the manifold $X=Z\# k\overline{\C P^{2}}\# l(S^1\times S^3)$ is given by
\begin{equation}\label{eq:cx}
\PP_X= -4\pi\sqrt{2c_1^2(Z)} \ .
\end{equation}
The supremum defining $\PP_X$ is realized by a metric if and only if $k=l=0$ and $Z$ admits a K\"ahler--Einstein metric 
of non-positive scalar curvature.
\end{thm}
We also give calculations and estimates for some other classes of four-manifolds, see in particular Theorems~\ref{t:symp},
\ref{t:kK3}, and~\ref{t:2}. Our results are applications, and, in some cases, generalizations of the recent results of F.~Fang 
and Y.~Zhang~\cite{FZ}, who discovered a relationship between Perelman's invariant and solutions of the Seiberg--Witten
monopole equations on four-manifolds. In particular, Fang and Zhang~\cite{FZ} already noted that Perelman's invariant is 
not a homeomorphism invariant. Our Theorems~\ref{t:kK3}, \ref{t:2} and~\ref{t3} contain quantitative results elaborating on 
this observation.

As a consequence of Seiberg--Witten theory it is now well known that the Yamabe invariant is in fact sensitive to the smooth
structure of a four-manifold. As was pointed out in~\cite{entropies}, this is also true for Gromov's minimal volume, although
this is harder to prove than for the Yamabe invariant. In Theorem~\ref{t:2} below we prove that the vanishing or non-vanishing of the minimal volume does depend on the smooth structure, and then discuss the relationship between this result and Perelman's invariant.

\newpage 

\section{Perelman's invariant}\label{s:perl}

We recall the basic definition from~\cite{P1}, compare also~\cite{KL,FZ}.

Let $M$ be a closed oriented manifold of dimension $n\geq 3$. For a Riemannian metric $g$ on $M$ and a function
$f\in C^{\infty}(M)$, Perelman defines
$$
\FF(g,f) = \int_M(s_g+\vert\nabla f\vert^2)e^{-f}\dvol_g \ ,
$$
where $s_g$ is the scalar curvature of $g$. Then Perelman's invariant of the Riemannian metric $g$ is
$$
\lambda_M(g)=\inf_{f\in C^{\infty}(M)}\{ \FF(g,f) \ \vert \ \int_M e^{-f}\dvol_g = 1\} \ .
$$
This infimum is actually a minimum, because it coincides with the smallest eigenvalue of the operator $d^*d+s_g$. 
It follows that $\lambda_M(g)$ depends continuously on $g$. If the scalar curvature $s_g$ is constant, then 
$\lambda_M(g)=s_g$.

The quantity 
$$
\PP_M(g)=\lambda_M(g)\cdot Vol(M,g)^{2/n}
$$
is scale invariant, and can be used to define a diffeomorphism invariant of $M$ by setting
\begin{equation}\label{eq:PP}
\PP_M = \sup_{g}\PP_M(g) \ .
\end{equation}
We shall call $\PP_M$ Perelman's invariant of $M$.

If the supremum in~\eqref{eq:PP} is achieved, then the corresponding metric has to be Einstein, cf.~\cite{P1,KL}.
We will see that very often this is not possible, so that the supremum is not a maximum.

Recall that for a conformal class $C$ on $M$, the Yamabe invariant of $C$ is defined by
$$
\mu_M(C) = \inf_{g\in C}\frac{\int_Ms_g\dvol_g}{Vol(M,g)^{(n-2)/n}} \ .
$$
By the solution to the Yamabe problem due to H.~Yamabe, N.~Trudinger, T.~Aubin and R.~Schoen, this infimum
is always realized by some metric $g_0$ of constant scalar curvature $s_0$. Such a metric is called
a Yamabe  metric. For a Yamabe metric $g_0\in C$ we  have
$$
\mu_M(C) = \frac{\int_Ms_0dvol_{g_0}}{Vol(M,g_0)^{(n-2)/n}} = s_0\cdot  Vol(M,g_0)^{2/n}=\PP_M(g_0) \ .
$$
The Yamabe invariant of $M$ is defined as  
$$
\mu_M = \sup_{C}\mu_M(C) \ .
$$
By the above calculation this can be written as 
$$
\mu_M = \sup \{ \PP_M(g_0) \ \vert \ g_0 \ \textrm{a Yamabe metric}\} \ .
$$
Comparing this with the definition~\eqref{eq:PP} of Perelman's invariant, we obtain the fundamental inequality
\begin{equation}\label{eq:YP}
\mu_M\leq\PP_M
\end{equation}
between the Yamabe and Perelman  invariants of $M$.

\section{Monopole classes}\label{s:mono}

Consider now a closed smooth oriented four-manifold $X$ with a $\Spc$-structure 
$\fs$. For every choice of Riemannian metric $g$, the Seiberg--Witten 
monopole equations for $(X,\fs)$ with respect to $g$ are a system of coupled 
equations for a pair $(A,\Phi)$, where $A$ is a $\Spc$-connection in the spin 
bundle for $\fs$ and $\Phi$ is a section of the positive spin bundle $V_{+}$.
The equations are:
\begin{equation}\label{eq:D}
     \DA \Phi =0 \ ,
\end{equation}
\begin{equation}\label{eq:curv}
     F_{\hat A}^+ =\sigma (\Phi ,\Phi ) \ ,
\end{equation}
with $\DA$ the half-Dirac operator defined on spinors of positive chirality, 
and ${\hat A}$ the connection induced by $A$ on the determinant of the spin 
bundle. The right-hand side of the curvature equation~\eqref{eq:curv} is the 
$2$-form which, under the Clifford module structure determined by $\fs$, 
corresponds to the trace-free part of $\Phi\otimes\Phi^{*}$.

Solutions $(A,\Phi)$ with $\Phi\equiv 0$ are called {\it reducible}. If 
there is a reducible solution, then $c=c_{1}(\fs)$ is represented by an 
anti-self-dual harmonic form because of~\eqref{eq:curv}. This implies 
$c^{2}\leq 0$, with equality if and only if $c$ is a torsion class.

The following definition is due to P.~Kronheimer~\cite{Kr}, see also~\cite{monopoles}.
\begin{defn}
    A class $c\in H^{2}(X,\Z)$ is called a {\it monopole class}, if there 
    is a $\Spc$-structure $\fs$ on $X$ with $c=c_{1}(\fs)$ for which the 
    monopole equations~\eqref{eq:D} and~\eqref{eq:curv} admit a solution 
    $(A,\Phi)$ for every choice of metric $g$. 
    \end{defn}
Of course, on manifolds for which the Seiberg--Witten invariants are 
well-defined, every basic class is a monopole class. The rationale for 
considering the concept of a monopole class is that the existence of 
solutions to the monopole equations has immediate consequences, even when 
the corresponding invariants vanish. 

The following result from~\cite{FZ} establishes an important relation between the 
Seiberg--Witten equations and Perelman's invariant:
\begin{prop}[Fang--Zhang~\cite{FZ}]\label{p:FZ}
If the monopole equations~\eqref{eq:D} and~\eqref{eq:curv} for $(X,\fs)$ with 
respect to $g$ admit an irreducible solution $(A,\Phi)$, then
\begin{equation}\label{eq:FZ}
\PP_X(g)\leq -4\pi\sqrt{2(c_1^+)^2} \ ,
\end{equation}
where $c_1^+$ denotes the projection of $c_1(\fs)$ to the $g$-self-dual subspace 
of $H^2(X;\R)$. 

If $c_1^+\neq 0$, then equality in~\eqref{eq:FZ} can only occur if $g$ is a 
K\"ahler metric of constant negative scalar curvature.
\end{prop}
The proof is an adaptation of the usual scalar curvature estimate for solutions of the 
Seiberg--Witten equations obtained by combining the Bochner--Weitzenb\"ock formula
with the equations, cf.~\cite{Wi}.

The following is a slight generalization of Theorem~1.1 in~\cite{FZ}:
\begin{thm}
Let $X$ be a smooth closed oriented $4$-manifold with a monopole class $c$  that is 
not a torsion class and satisfies $c^2\geq 0$. Then 
\begin{equation}\label{eq:bound}
\PP_X\leq -4\pi\sqrt{2c^2} \ .
\end{equation}
\end{thm}
\begin{proof}
Because $c$ is assumed to be a monopole class, there is a $\Spc$-structure $\fs$ with
$c_1(\fs)=c$ such that the monopole equations for $\fs$ have a solution $(A,\Phi)$ for 
every choice of metric $g$. If the solution is irreducible, i.~e. $\Phi\neq 0$, then~\eqref{eq:FZ}
implies 
$$
\PP_X(g)\leq -4\pi\sqrt{2(c^+)^2} \leq -4\pi\sqrt{2c^2} \ .
$$
If $c^2>0$, then it is clear that all solutions are irreducible. If $c^2=0$, then we use the 
assumption that $c$ is not a torsion class, to argue that solutions must be irreducible for generic $g$. In fact,
this follows from:
\begin{lem}[Donaldson]\label{l:D}
    Let $c\in H^{2}(X,\Z)$ be a non-torsion class. If $b_{2}^{+}(X)>0$, then for a 
    generic metric $g$, there is no anti-self-dual harmonic form 
    representing the image of $c$ in $H^{2}(X,\R)$.
    \end{lem}
A proof of the lemma can be found in~\cite{DK}. In our case the intersection form must be 
indefinite because there is a non-torsion class of square zero, and thus $b_{2}^{+}(X)>0$ 
holds. Therefore we  obtain the desired estimate for generic metrics $g$. As Perelman's 
invariant $\PP_X(g)$ depends continuously on $g$, the estimate holds for all $g$.
\end{proof}

Next we adapt the proof of Theorem~4.2 in~\cite{monopoles} to show that in the presence of a 
monopole class, Perelman's invariant can be used 
to bound the number of smooth exceptional spheres in a four-manifold.
\begin{thm}\label{t:split}
    Let $X$ be a smooth four-manifold with a monopole class $c$. The maximal number $k$ of 
    copies of $\overline{\C P^{2}}$ that can be split off smoothly is bounded by
    \begin{equation}\label{eq:split}
    k\leq\frac{1}{32\pi^2}\PP_X^2-c^2 \ .
    \end{equation}
    \end{thm}
\begin{proof}
    Suppose that $X\cong Y\#k\overline{\C P^{2}}$, with $k>0$, and write 
    $c=c_{Y}+\sum_{i=1}^{k}a_{i}e_{i}$, with respect to the obvious 
    direct sum decomposition of $H^{2}(X,\Z)$. Here $e_{i}$ are the 
    generators for the cohomology of the $\overline{\C P^{2}}$ summands. 
    Note that the $a_{i}$ are odd integers because $c$ must be 
    characteristic. In particular $c$ can not be a torsion class.     
    
    Now the reflections in the $e_{i}$ are realised by 
    self-diffeomorphisms of $X$, and the images of our monopole class under 
    these diffeomorphisms are again monopole classes. Thus, moving $c$ 
    by a diffeomorphism, we can arbitrarily change $e_{i}$ to its negative.
    
    Given a metric $g$ on $X$, we choose the signs in such a 
    way that $a_{i}e_{i}^{+}\cdot c_{Y}^{+}\geq 0$. Then we find
$$
    (c^{+})^{2} = \left( c_{Y}^{+}+\sum_{i=1}^{k}a_{i}e_{i}^{+}\right)^{2}\geq 
    (c_{Y}^{+})^{2}\\
    \geq c_{Y}^{2}=c^{2}+\sum_{i=1}^{k}a_{i}^{2}  \ .
$$
If $g$ is generic, then there are irreducible solutions to the monopole equations,
and applying~\eqref{eq:FZ} to $c$ and $g$ gives
    $$
    (c^{+})^{2}\leq\frac{1}{32\pi^2}(\PP_X(g))^2 \ .
    $$
    Combining the two inequalities and noting that $a_{i}^{2}\geq 1$ 
    because all the $a_{i}$ are odd integers shows
$$
        k\leq\frac{1}{32\pi^2}(\PP_X(g))^2-c^2 \ 
$$
for generic $g$. By continuity of Perelman's invariant this holds for all $g$, giving~\eqref{eq:split}.
    \end{proof}
It may not be obvious why this Theorem is interesting, but this should become clear
by looking at the following special case:
\begin{thm}\label{t:symp}
Let $Z$ be a minimal symplectic four-manifold with $b_2^+(Z)>1$, and $X=Z\# k\overline{\C P^{2}}\# l(S^1\times S^3)$.
Then 
\begin{equation}\label{eq:symp}
\PP_X\leq -4\pi\sqrt{2c_1^2(Z)} \ .
\end{equation}
\end{thm}
The case $l=0$ was previously proved by Fang and Zhang in Theorem~1.4 of~\cite{FZ}.
\begin{proof}
Let $Y=Z\# k\overline{\C P^{2}}$, which we  can think of as a symplectic blowup. By the result of 
C.~Taubes~\cite{T}, the first Chern class $c_1(Y)$ of a symplectic structure is a Seiberg--Witten basic class with
numerical Seiberg--Witten invariant $\pm 1$. In particular, it is a monopole class on $Y$.

	Now consider $X=Y\# l(S^{1}\times S^{3})$. Although its numerical Seiberg--Witten invariants 
	must vanish, cf.~\cite{KMT,Bourbaki}, we claim that each 
	of the basic classes with numerical Seiberg--Witten invariant $=\pm 1$ 
	on $Y$ gives rise to a monopole class on $Z$. 
	There are two ways to see this. One can extract our claim from the 
	connected sum formula~\cite{B} for the stable cohomotopy refinement of 
	Seiberg--Witten invariants introduced by S.~Bauer and M.~Furuta~\cite{BF}, 
	cf.~\cite{F}. Alternatively, one uses the invariant defined by 
	the homology class of the moduli space of solutions to the monopole 
	equations, as in~\cite{Bourbaki}. This means that the first homology 
	of the manifold is used, and here this is enough to obtain a non-vanishing 
	invariant. Using this invariant, our claim follows from Proposition~2.2 
	of P.~Ozsv\'ath and Z.~Szab\'o~\cite{OS}. See also K.~Froyshov~\cite{Fro}.

We apply~\eqref{eq:split} to $X$ with the monopole class $c_1(Y)\in H^2(Y;\Z)=H^2(X;\Z)$ to 
obtain 
$$
k\leq \frac{1}{32\pi^2}\PP_X^2-c_1^2(Y) = \frac{1}{32\pi^2}\PP_X^2-c_1^2(Z)+k \ . 
$$
As $\PP_X$ is non-positive in this case, the claim follows.
\end{proof}
In this argument the expected dimension of the moduli space of solutions to the Seiberg--Witten
equations is $=0$ on $Y$, but is $=l$ on $X$. As explained in~\cite{monopoles}, 
results like Theorem~\ref{t:split} are stronger the larger the expected dimension of the 
moduli space is.

In the case of K\"ahlerian complex surfaces we can combine this upper bound for Perelman's invariant 
with the lower bound given by the Yamabe invariant to obtain a complete calculation.
\begin{proof}[Proof of Theorem~\ref{t:cx}]
Every compact complex surface $Z$ with even first Betti number admits a K\"ahler structure, and is therefore 
symplectic. For simplicity we are assuming $b_2^+(Z)>1$, so that surfaces of negative Kodaira dimension
do not occur. For $Z$ with a K\"ahler structure of non-negative Kodaira dimension, holomorphic and symplectic
minimality coincide, cf.~\cite{HK}. Thus $Z$ is symplectically minimal and we can apply Theorem~\ref{t:symp} 
to obtain 
\begin{equation}\label{eq:inter}
\PP_X\leq -4\pi\sqrt{2c_1^2(Z)} \ .
\end{equation}

For the reverse inequality consider first the case $k=l=0$, i.~e. $X=Z$. Then we have $\PP_Z\geq\mu_Z$ 
by~\eqref{eq:YP}, and, if $Z$ is of general type, then the Yamabe invariant $\mu_Z$ equals $-4\pi\sqrt{2c_1^2(Z)}$ 
by the result of C.~LeBrun~\cite{lebrun}. This is easy to see when $Z$ admits a K\"ahler--Einstein metric of negative 
scalar curvature, cf.~\cite{lebrun,FZ}. In the case when such a metric does not exist, one has to consider sequences 
of metrics which suitably approximate an orbifold K\"ahler--Einstein metric on the canonical model of $Z$, see~\cite{lebrun}. 
If $Z$ is not of general type, then $c_1^2(Z)=0$, and $Z$ is either a $K3$ surface, an elliptic surface, or an Abelian surface, 
cf.~\cite{BPV}. In all these cases $Z$ does not admit a metric of positive scalar curvature, but it does collapse with bounded 
scalar curvature, in fact even with bounded Ricci curvature, see~LeBrun~\cite{lebrun3}. Collapsing with bounded scalar 
curvature can also be seen from the result  of G.~Paternain and J.~Petean~\cite{PP} that $Z$ has an $\FF$-structure. 
Because $Z$ collapses with bounded scalar curvature, its Yamabe and Perelman invariants vanish. This completes 
the proof of~\eqref{eq:cx} in the case $k=l=0$. 

Now we allow $k$ and $l$ to be positive. By Proposition~4.1 of Fang and Zhang~\cite{FZ}, Perelman's invariant does 
not decrease under connected sum with $\overline{\C P^{2}}$ and with $S^1\times S^3$. As the upper bound~\eqref{eq:inter}
is achieved on $Z$ and is unchanged by the connected sum, we conclude that it is an equality for all positive $k$ and $l$
as well.

Finally we discuss the question whether the supremum in the definition of $\PP_X$ is a maximum. If this is the case,
then the corresponding metric on $X$ is an Einstein metric. If $Z$ is of general type, then $c_1^2(Z)>0$, and we can use 
the discussion of the equality case in Proposition~\ref{p:FZ} to conclude that the critical metric on $X$ is K\"ahler as well, 
and the scalar curvature is negative. This implies $k=l=0$. If $Z$ is not of general type then $c_1^2(Z)=0$. Now the 
Hitchin--Thorpe inequality~\cite{HT} for the Einstein metric implies $k=l=0$. We are then in the case of equality of the 
Hitchin--Thorpe inequality, and $X=Z$ is Ricci-flat K\"ahler. 
\end{proof}
\begin{rem}
Instead of using the behaviour of Perelman's invariant under connected sum, one can alternatively argue with 
the corresponding result for the Yamabe invariant (and for $\FF$-structures).
\end{rem}
\begin{rem}
Theorems~\ref{t:symp} and~\ref{t:cx} have extensions to the case of manifolds with $b_2^+=1$.
For the latter one also has to consider rational and ruled symplectic manifolds, which do admit metrics of positive scalar
curvature and therefore have positive Perelman invariant.
\end{rem}

\section{Examples and applications}\label{s:ex}

In this section we give some examples illustrating the estimates and calculations of Perelman's invariant, with 
special emphasis on its dependence on the smooth structure.

First, we have the following:
\begin{thm}
The number of distinct values that Perelman's invariant can take on the smooth structures in a fixed homeomorphism type of simply 
connected four-manifolds is unbounded.
\end{thm}
\begin{proof}
By the standard geography results for minimal surfaces of general type going back to U.~Persson~\cite{per}, compare also~\cite{BPV},
we can do the following: for every positive integer $n$ we find positive integers $x$ and $y$ with the properties that all pairs
$(x-i,y+i)$ with $i$ running from $1$ to $n$ are realized as $(c_2(Z_i),c_1^2(Z_i))$ for some simply connected minimal complex
surface $Z_i$ of general type. Let $X_i$ be the $i$-fold blowup of $Z_i$. Then all the $X_i$ for $i$ from $1$ to $n$ are 
simply connected and non-spin and have the same Chern numbers $(x,y)$. Therefore, by M.~Freedman's result~\cite{freed} they are  
homeomorphic to each other. However, by the above Theorem~\ref{t:cx}, the $X_i$ have pairwise different Perelman invariants.
\end{proof}
While this construction does produce arbitrarily large numbers of distinct Perelman invariants among homeomorphic four-manifolds,
it can never produce infinitely many. Of course it is now easy to construct manifolds with infinitely many distinct smooth structures
admitting symplectic forms. However, Theorem~\ref{t:symp} does not seem to be strong enough to show that their Perelman 
invariants take on infinitely many values. Therefore, we leave open the following:
\begin{conj}\label{con}
On a suitable homeomorphism type, the Perelman invariant takes on infinitely many distinct values.
\end{conj}

There are also spin manifolds with smooth structures with several distinct Perelman invariants:
\begin{thm}\label{t:kK3}
The manifold $X=3K3\# 4(S^2\times S^2)$ has:
\begin{itemize}
\item a smooth structure $X_0$ with $\PP_{X_0}=0$,
\item a smooth structures $X_1$ with $\PP_{X_1}=-16\pi$,
\item infinitely many smooth structures $X_i$ with $\PP_{X_i}\leq-16\sqrt2\pi<-16\pi$.
\end{itemize}
The supremum defining the Perelman invariant is attained for $X_1$, but not for $X_0$ and the $X_i$.
\end{thm}
\begin{proof}
The smooth structure $X_0$ is the standard one given by the connected sum. By the Lichnerowicz argument it has no 
metric of positive scalar curvature. As it does collapse with bounded scalar curvature, we conclude $\PP_{X_0}=0$.
If this supremum were attained, then the corresponding metric would have to be Ricci-flat. As the signature is non-zero,
we would have a parallel harmonic spinor, showing that the manifold is K\"ahler, which is clearly not possible\footnote{See the 
appendix to~\cite{monopoles} for details of this argument.}.

The smooth structure $X_1$ underlies the complex algebraic surface obtained as the double cover of the projective 
plane branched in a smooth holomorphic curve of degree $10$. This is homeomorphic to $X$ by~\cite{freed}. 
By Theorem~\ref{t:cx}, $\PP_{X_1}=-16\pi$. As the canonical bundle of $X_1$ is ample, $X_1$ has a K\"ahler--Einstein metric
of negative scalar curvature by the results of T.~Aubin~\cite{A} and S.-T.~Yau~\cite{Y} on the Calabi conjecture. This metric
achieves the supremum for the Perelman invariant.

The smooth structures $X_i$ are constructed as follows, cf.~\cite{monopoles}.
Let $M$ be a symplectic spin manifold with $c_{1}^{2}(M)=16$ and $\chi(M)=4$, where $\chi=\frac{1}{4}(e+\sigma)$ 
denotes the holomorphic Euler characteristic. 
    Such manifolds exist by the results of D.~Park and Z.~Szab\'o~\cite{PS}. By 
    Freedman's classification~\cite{freed}, such an $M$ is homeomorphic to 
    $K3\# 4(S^{2}\times S^{2})$. Take $N=K3$, and $O$ 
    the symplectic spin manifold obtained from $K3$ by performing a logarithmic 
    transformation of odd multiplicity $i$. 
    By the connected sum formula for the stable cohomotopy refinement of the Seiberg--Witten
   invariants due to Bauer and Furuta~\cite{BF,B}, the connected sum $M\# N\# O$ has monopole classes $c$
   which are the sums of the basic classes on the different summands.  Note that $c^2=16$. Therefore 
    $\PP_{X_i}\leq-16\sqrt2\pi$ by Theorem~\ref{t:symp}. It was shown in~\cite{monopoles} that 
    as we increase $i$, the multiplicity of the logarithmic transformation, we do indeed get infinitely 
    many distinct smooth structures. It was also shown in~\cite{monopoles} that the $X_i$ do not admit 
    any Einstein metrics. Therefore, the supremum for the Perelman invariant can not be achieved for them.
    \end{proof}
\begin{rem}
    The manifold $X$ has another infinite sequence of smooth structures, 
    which are distinct from the ones discussed above. 
    R.~Fintushel and R.~Stern~\cite{FScusp} have shown 
    that one can perform cusp surgery on a torus in $S$ to construct 
    infinitely many distinct smooth structures which are irreducible and 
    non-complex, and are therefore distinct from the smooth structures we 
    consider. These smooth structures have negative Perelman invariants, and it is 
    unknown whether they admit Einstein metrics.
\end{rem}

One can give many similar examples on larger manifolds.
We can even obtain interesting results for parallelizable manifolds:
\begin{thm}\label{t:2}
    For every $k\geq 0$ the manifold $X_{k}=k(S^{2}\times S^{2})\# (1+k)(S^{1}\times 
	S^{3})$ with its standard smooth structure has zero minimal volume and Perelman
	invariant $\PP_{X_k}=+\infty$.
	
	If $k$ is odd and large enough, then there are infinitely many pairwise 
	non-diffeomorphic smooth manifolds $Y_{k}$ homeomorphic to $X_{k}$, all of 
	which have strictly positive minimal volume and strictly negative Perelman invariant.
	Moreover, the supremum in the definition of the Perelman invariant is not achieved.
   All the $Y_k$ have the property that $Y_k\# (S^2\times S^2)$ is diffeomorphic to 
	$X_k\# (S^2\times S^2)$.
    \end{thm}
    \begin{proof}
	Note that $X_{0}=S^{1}\times S^{3}$ has obvious free circle actions, and 
	therefore collapses with bounded sectional curvature. To see that all 
	$X_{k}$ have vanishing minimal volume it suffices to construct 
	fixed-point-free circle actions on them, cf.~M.~Gromov~\cite{gromov}.
	
	The product $S^{2}\times S^{2}$ has a diagonal effective circle action 
	which on each factor is rotation around the north-south axis. It has 
	four fixed pints, and the linearization of the action induces one 
	orientation at two of the fixed points, and the other orientation at 
	the remaining two. The induced action on the boundary of an $S^{1}$-invariant 
	small ball around each of the fixed points is the Hopf action on $S^{3}$. 
	By taking equivariant connected sums at fixed points, pairing fixed points at 
	which the linearizations give opposite orientations, we obtain effective circle 
	actions with $2+2k$ fixed points on the connected sum $k(S^{2}\times S^{2})$ 
	for every $k\geq 1$. Now we have $1+k$ fixed points at which the 
	linearization induces one orientation, and $1+k$ at which it induces 
	the other orientation. Then making equivariant self-connected sums at pairs 
	of fixed points with linearizations inducing opposite orientations we 
	finally obtain a free circle action on 
	$X_{k}=k(S^{2}\times S^{2})\# (1+k)(S^{1}\times S^{3})$.
	
	That the Perelman invariant of $X_k$ is infinite follows from the fact that this is 
	so for $S^2\times S^2$.

	If $k$ is odd and large enough, then there are 
	symplectic manifolds $Z_{k}$ homeomorphic (but not 
	diffeomorphic) to $k(S^{2}\times S^{2})$, see for example B.~Hanke, J.~Wehrheim and myself~\cite{HKW}. 
	By the construction given in~\cite{HKW}, we may assume that $Z_{k}$ 
	contains the Gompf nucleus of an elliptic surface. By performing 
	logarithmic transformations inside this nucleus, we can vary the 
	smooth structures on the $Z_{k}$ in such a way that the number of 
	Seiberg--Witten basic classes with numerical Seiberg--Witten 
	invariant $=\pm 1$ becomes arbitrarily large, cf.~Theorem~8.7 
	of~\cite{FSrat} and Example~3.5 of~\cite{monopoles}.
	
	Consider $Y_{k}=Z_{k}\# (1+k)(S^{1}\times S^{3})$. This is clearly 
	homeomorphic to $X_{k}$. The basic classes on $Z_k$ give rise
	to monopole classes which are special in the sense of~\cite{monopoles}.
	As their number is unbounded, we have infinitely many distinct smooth 
	structures. By the construction given in~\cite{HKW}, the $Z_k$ dissolve 
	after a single stabilization with $S^2\times S^2$, and therefore the same is 
	true for the $Y_k$.

	As $Y_{k}$ has non-torsion monopole classes $c$ with 
	$c^{2}=2\chi(Z_{k})+3\sigma(Z_{k})=4+4k>0$, we find that $\PP_{Y_k}\leq -8\pi\sqrt{2(1+k)}$
	from Theorem~\ref{t:symp}. This supremum is not achieved, because $Y_k$ is parallelizable 
	but not flat, and can therefore not carry an Einstein metric by the discussion of the equality  
	case of the Hitchin--Thorpe inequality~\cite{HT}.
	Note that $Y_{k}$ can not collapse with bounded scalar curvature. 
	{\it A fortiori} it cannot collapse with bounded 
	sectional curvature, and so its minimal volume is strictly positive. 
	\end{proof}

\begin{rem}
That connected sums of manifolds with vanishing minimal volumes may have non-vanishing
minimal volumes is immediate by looking at connected sums of tori.
The manifolds $X_k$ discussed above have the property that their minimal volumes
vanish, although they are connected sums of manifolds with non-vanishing minimal volumes.
Thus the minimal volume, and even its (non-)vanishing, does not behave in a straightforward
manner under connected sums.
\end{rem}
This remark was motivated by the recent preprint of G.~Paternain and J.~Petean~\cite{PPnew}.
After an earlier version of Theorem~\ref{t:2} appeared on the arXiv in~\cite{entropies}, 
these authors remarked on the complicated behaviour of the minimal volume under connected sums
based on some $6$-dimensional examples, see Remark~3.1 in~\cite{PPnew}.

\begin{rem}
The insistance on manifolds that dissolve after a single stabilization with $S^2\times S^2$ was originally 
motivated by C.~T.~C.~Wall's theorem~\cite{W} showing that every simply connected smooth four-manifold 
dissolves after some number of stabilizations with $S^2\times S^2$. On the one hand, because of work of 
R.~Mandelbaum, B.~Moishezon and R.~Gompf, and also because of~\cite{BK,HKW}, we now know that one
stabilization often suffices. On the other hand, there are no examples where it is known that one stabilization
does not suffice. Perelman's invariant sheds some unexpected light on this, because after a single stabilization
with $S^2\times S^2$ Perelman's invariant tends to become infinite. In particular, after a single stabilization there 
are no more non-torsion monopole classes. It is tempting to speculate that the Ricci flow might be useful in
investigating the question whether manifolds do indeed always dissolve after a single stabilization.
\end{rem}

Going in the opposite direction of Conjecture~\ref{con}, we have the following:
\begin{thm}\label{t3}
There are homeomorphism types of simply connected four-manifolds that contain infinitely many distinct smooth structures
with the same Perelman invariant.
\end{thm}
\begin{proof}
The easiest example is furnished by the homeomorphism type underlying the $K3$ surface. This is spin with non-zero signature,
and so by the Lichnerowicz argument no manifold in this homeomorphism type admits a metric of positive scalar curvature. 
However, there are infinitely many smooth structures underlying complex elliptic surfaces in this homeomorphism type. 
As all these collapse with bounded scalar curvature, their Perelman invariants vanish.
\end{proof}
In this case the supremum for the Perelman invariant is attained for the standard smooth structure, but not for any 
of the other ones. This is because the standard smooth structure is the only one admitting an Einstein metric. These
were in fact the first examples showing that the existence of an Einstein metric depends on the smooth structure, cf.~\cite{K}.
\begin{rem}
If we only look for arbitrarily large  numbers of distinct smooth structures, rather than for infinitely many, then we can choose examples
with non-vanishing Perelman invariants. For example, V.~Braungardt and I proved in~\cite{BK} that there are arbitrarily
large tuples of non-diffeomorphic minimal surfaces of general type with ample canonical bundles. By Theorem~\ref{t:cx}
above on any such tuple the Perelman invariant is a negative constant. As shown in~\cite{BK}, these examples can be  chosen
to be spin or non-spin. In the non-spin case they also have infinitely many other, non-complex, smooth structures
which, by Theorem~\ref{t:symp}, have even smaller Perelman invariant.
\end{rem}


\bibliographystyle{amsplain}

\bigskip

\end{document}